\documentclass[12pt]{article}
\usepackage{latexsym}
\usepackage[francais]{babel}
\usepackage{mathrsfs}
\usepackage{a4wide}
\usepackage{theorem}
\usepackage{amsfonts}
\usepackage{amsmath}
\usepackage{amscd}
\usepackage{times}

\theorembodyfont{\sl}
\newtheorem{theorem}{Th\'eor\`eme}

\theorembodyfont{\rmfamily}
\newtheorem{rem}{Remarque}

\newenvironment{prf}[1]{\trivlist
\item[\hskip \labelsep{\it #1}]}{~\hspace{\fill}~$\square$\endtrivlist}
\newenvironment{proof}{\begin{prf}{\bf Preuve}}{\end{prf}}

\newcommand{\square}{\Box}
\newcommand{\ZZ}{{\mathbb Z}}
\newcommand{\FF}{{\mathbb F}}
\newcommand{\PP}{{\mathbb P}}
\newcommand{\GL}{{\rm GL}}
\newcommand{\SL}{{\rm SL}}
\newcommand{\CC}{{\mathbb C}}
\newcommand{\QQ}{{\mathbb Q}}
\newcommand{\Qbar}{{\overline{\QQ}}}

\newcommand{\ld}{\langle}
\newcommand{\rd}{\rangle}
\newcommand{\Aut}{{\rm Aut}}
\newcommand{\Hom}{{\rm Hom}}

\newcommand{\nonsplit}{{\mbox{\rm\scriptsize non-split}}}
\newcommand{\jac}{{\rm jac}}

\newcommand{\ol}{\overline}
\newcommand{\End}{{\rm End}}

\newcommand{\calC}{{\mathscr C}}
\newcommand{\calA}{{\mathscr A}}

\newcommand{\FFq}{{\FF\,'}}
\newcommand{\cprod}{\sqcup}
\newcommand{\Sym}{{\rm Sym}}
\newcommand{\ZZl}{\ZZ_{(l)}}

\begin{document}
\title{Sur un r\'esultat d'Imin Chen}
\author{Bart de Smit\footnote{
boursier du Koninklijke Nederlandse Akademie van Wetenschappen.}\and
Bas Edixhoven\footnote{membre de l'Institut Universitaire
de France, b\'en\'eficiaire du European TMR Network
Contract ERB FMRX 960006 ``arithmetic algebraic geometry''.}}
\date{9 f\'evrier 2000}

\maketitle

\begin{abstract}
Nous g\'en\'eralisons un r\'esultat de Chen concernant une isog\'enie
entre produits de jacobiennes de courbes modulaires associ\'ees \`a
des sous-groupes de~$\GL_2(\FF_p)$. Cette g\'en\'erali\-sation
s'applique \`a des objets munis d'une action de $\GL_2(\FF)$, avec
$\FF$ un corps fini quelconque, dans des cat\'egories additives plus
g\'en\'erales.
\end{abstract}

\def\abstractname{R\'esum\'e (version anglaise)}%
\begin{abstract}
We generalize a result of Chen concerning an isogeny between products
of jacobians of modular curves associated to subgroups
of~$\GL_2(\FF_p)$. This generalization concerns objects with an action
by $\GL_2(\FF)$, with $\FF$ an arbitrary finite field, in more general
additive categories.
\end{abstract}

\bigskip\noindent Le but de cet article est de red\'emontrer et de
g\'en\'eraliser un r\'esultat d'Imin Chen, concernant certaines
identit\'es entre fonctions z\^eta de courbes modulaires, ou, ce qui
revient au m\^eme, certaines isog\'enies entre des produits de
jacobiennes de telles courbes.

Soient $\FF$ un corps fini, $\FFq$ une extension quadratique de $\FF$
et $G:=\GL_2(\FF)$.  Le groupe $G$ agit sur $\PP^1(\FF)$ et sur
$\PP^1(\FFq)$.  Soit $B$ le fixateur dans $G$ d'un point $\infty$ dans
$\PP^1(\FF)$, $T$ le fixateur dans $G$ de deux points $0,
\infty\in\PP^1(\FF)$, $N$ le fixateur de l'ensemble $\{0, \infty\}$,
$T'$ le fixateur d'un point dans $\PP^1(\FFq)-\PP^1(\FF)$, et $N'$ le
fixateur de sa $\Aut_\FF(\FFq)$-orbite.  On appelle $B$ un sous-groupe
de Borel, $T$ un tore maximal d\'eploy\'e, et $T'$ un tore maximal non
d\'eploy\'e.  Le groupe $N'$ est le normalisateur de $T'$ et si
$\FF\ne \FF_2$, $N$ est le normalisateur de $T$.

Rappelons d'abord le r\'esultat de Chen. Pour $n\geq1$ entier, soit
$X(n)_\QQ$ la courbe modulaire qui est l'espace de modules (grossier
si $n<3$) compactifi\'e de paires $(E/S/\QQ,\phi)$, o\`u $S$ est un
$\QQ$-sch\'ema, $E/S$ une courbe elliptique et $\phi\colon
(\ZZ/n\ZZ)^2\to E[n]$ un isomorphisme de $S$-sch\'emas en groupes
(voir \cite{Katz-Mazur}, o\`u ce probl\`eme de modules est appel\'e
``naive level $n$ structure''). Par construction, le groupe
$\GL_2(\ZZ/n\ZZ)$ agit \`a droite sur $X(n)_\QQ$~: un \'el\'ement $g$
envoie $(E/S/\QQ,\phi)$ vers $(E/S/\QQ,\phi\circ g)$. Cette action
induit, par fonctorialit\'e de Picard, une action \`a gauche sur la
jacobienne de~$X(n)_\QQ$.  La jacobienne d'une courbe $X$ sur $\QQ$
est un sch\'ema ab\'elien sur $\QQ$ not\'ee $\jac(X)$.  Le r\'esultat
de Chen est alors le suivant (voir \cite[Theorem~1 et \S10]{Chen}).

\begin{theorem}[Chen, 1994]\label{thm1.5}
Si $\FF=\FF_p$ avec $p$ premier, $G$ agit sur $X=X(p)_\QQ$.  Alors,
$\jac(X/T)$ est isog\`ene \`a $\jac(X/T')\times\jac(X/B)^2$ et
$\jac(X/N)$ est isog\`ene \`a $\jac(X/N')\times\jac(X/B)$.
\end{theorem}

\noindent
Certains quotients de $X$ sont connus sous d'autres noms (voir
\cite[Ch.~11]{Katz-Mazur}, \cite[p. 36]{Maz}):
$$
X/T\cong X_0(p^2)_\QQ, \quad
X/N'\cong X(p)_\nonsplit, \quad
X/B\cong X_0(p)_\QQ, \quad
X/N\cong X_0(p^2)_\QQ/\ld w_{p^2}\rd.
$$
La preuve donn\'ee par Chen consiste \`a montrer que les traces des
op\'erateurs de Hecke $T_n$ avec $n$ premier \`a $p$ sur les
jacobiennes en question satisfont les identit\'es requises pour
conclure, par la relation d'Eichler et Shimura et par la conjecture de
Tate (d\'emontr\'ee par Faltings), \`a l'existence d'une telle
isog\'enie.  Une preuve n'utilisant que la th\'eorie des
repr\'esentations de $G$ a \'et\'e donn\'ee dans~\cite{Edixhoven1};
cette note am\'eliore la m\'ethode de \cite{Edixhoven1}, et donne
quelques r\'esultats suppl\'ementaires (Th\'eor\`emes 3 et~4).  Voir
\cite{Dar} pour des applications aux variantes de l'\'equation de
Fermat.  Dans \cite{ChambertSauvageot} on trouve une preuve,
inspir\'ee par \cite{Edixhoven1} et par l'interpr\'etation ad\'elique
des formes modulaires.

Pour un groupe $H$ et $M$ un objet muni d'une action de $H$ dans une
cat\'egorie $\calC$ nous dirons que $M$ admet des invariants par $H$
si le foncteur $\Hom(-,M)^H$ est repr\'esentable. Dans ce cas, nous
notons $M^H$ le sous-objet de $M$ repr\'esentant ce foncteur.

\begin{theorem}\label{thm2.1}
Soit $\calC$ une cat\'egorie additive et $\QQ$-lin\'eaire.
Soit $M$ un objet de $\calC$ muni d'une action $\alpha$
de $G$, qui admet des invariants par les sous-groupes de~$G$.
Alors il existe des isomorphismes fonctoriels en~$(M,\alpha)$:
\begin{eqnarray}
M^T \oplus M^G \oplus M^G & \cong & M^{T'} \oplus M^B \oplus M^B \\
M^N \oplus M^G & \cong & M^{N'} \oplus M^B.
\end{eqnarray}
\end{theorem}

\noindent
Pour voir que le Th\'eor\`eme~\ref{thm1.5} est un cas sp\'ecial, on
applique le Th\'eor\`eme~\ref{thm2.1} avec $\calC=\QQ\otimes\calA$,
o\`u $\calA$ est la cat\'egorie des vari\'et\'es ab\'eliennes
sur~$\QQ$.  Les objets de $\calC$ sont donc les m\^emes que ceux de
$\calA$, et on a, pour deux objets $A$ et $B$,
$\Hom_{\calC}(A,B)=\QQ\otimes\Hom_{\calA}(A,B)$.  Pour un objet $A$ de
$\calA$ nous noterons $\QQ\otimes A$ l'objet correspondant
de~$\calC$. Par construction, $A$ et $B$ dans $\calA$ sont isog\`enes
si et seulement si $\QQ\otimes A$ et $\QQ\otimes B$ sont isomorphes
dans $\calC$. La cat\'egorie $\cal A$ est ab\'elienne (donc les objets
des invariants existent), $\QQ$-lin\'eaire et m\^eme semi-simple.
Soit maintenant $M:=\QQ\otimes\jac(X)$.  Pour tout sous-groupe $H$ de
$G$ on a alors $M^H=\QQ\otimes\jac(X/H)$.  On a $M^G=0$ car $X/G$ est
de genre z\'ero.

\begin{proof}
{\bf du Th\'eor\`eme 2} Par le lemme de Yoneda (voir \cite[Ch.~0,
\S1]{Grothendieck}), on se ram\`ene au cas o\`u $\calC$ est la
cat\'egorie des $\QQ$-espaces vectoriels. Notons, pour $H\subset G$ un
sous-groupe, $\QQ[G/H]$ le $\QQ[G]$-module donn\'e par l'action de
multiplication \`a gauche de $G$ sur~$G/H$~; c'est l'induite \`a $G$
de la repr\'esentation triviale de~$H$. Avec ces notations, on a, pour
tout $\QQ[G]$-module $M$~:
$$
M^H = \Hom_{\QQ[H]}(\QQ,M) =
\Hom_{\QQ[G]}(\QQ[G/H],M).
$$
Il suffit donc de montrer qu'il existe des isomorphismes
de $\QQ[G]$-modules~:
\begin{eqnarray}
\QQ[G/T] \oplus \QQ \oplus \QQ & \cong &
\QQ[G/T'] \oplus \QQ[G/B] \oplus \QQ[G/B]; \\
\QQ[G/N] \oplus \QQ & \cong &
\QQ[G/N'] \oplus \QQ[G/B].
\end{eqnarray}
D'apr\`es \cite[\S12]{Serre1}, il suffit, pour montrer que deux
$\QQ[G]$-modules de type fini sont isomorphes, de v\'erifier que leurs
caract\`eres co\"\i ncident. Dans notre cas, tout \'el\'ement de $G$ a
un conjugu\'e dans $B$ ou dans~$T'$. Pour montrer l'existence du
premier isomorphisme, nous allons montrer que les $G$-ensembles
$X:=G/T\cprod\{\cdot\}\cprod\{\cdot\}$ et $Y:=G/T'\cprod G/B\cprod
G/B$ sont isomorphes en tant que $B$-ensembles et en tant que
$T'$-ensembles.  Par la d\'efinition de $T$, $T'$ et $B$ on a:
$$
X=(\PP^1(\FF)\times\PP^1(\FF)-\Delta)
\cprod\{\cdot\}\cprod\{\cdot\}, \qquad
Y=\PP^1(\FFq)\cprod \PP^1(\FF).
$$
Notons $\sigma$ l'\'el\'ement non trivial de $\Aut_\FF(\FFq)$.  En
tant que $B$-ensemble, $\PP^1(\FF)$ est la r\'eunion disjointe de la
droite affine $\FF$ et d'un point~$\infty$. Comme l'action de
$B/\FF^*$ sur $\FF$ est simplement doublement transitive, on a un
isomorphisme de $B$-ensembles:
$$
X\cong B/\FF^* \cprod \FF \cprod \FF\cprod \{\cdot\}\cprod \{\cdot\}.
$$
D'autre part, $\PP^1(\FFq)-\PP^1(\FF)$ est un $B/\FF^*$-ensemble
libre, car si un \'el\'ement de $B$ en fixe $P$, il fixe les trois
points $P$, $\sigma(P)$ et $\infty$ de~$\PP^1(\FFq)$, donc il est
scalaire.  Ceci ach\`eve la preuve que $X$ et $Y$ sont isomorphes en
tant que $B$-ensembles. Consid\'erons maintenant $X$ et $Y$ comme
$T'$-ensembles.  Par d\'efinition, $T'$ fixe deux points conjugu\'es
de $\PP^1(\FFq)$, donc tout \'el\'ement de $T'$ qui fixe un
troisi\`eme point de $\PP^1(\FFq)$ est scalaire. Il en r\'esulte que
$X$ et $Y$ ont chacun deux points fixes par $T'$, dont les
compl\'ements sont des $T'/\FF^*$-ensembles libres. Comme $\# X=\# Y$,
on a $X\cong Y$ en tant que $T'$-ensembles, ce qui termine la preuve
de l'existence d'un isomorphisme~(3).

Pour \'etablir (4), il suffit de montrer que les $G$-ensembles:
$$
\ol{X}:= (\Sym^2(\PP^1(\FF))-\Delta) \cprod\{\cdot\}, \qquad
\ol{Y}:= \PP^1(\FFq)/\ld\sigma\rd
$$
sont isomorphes en tant que $B$-ensembles et $T'$-ensembles.  Comme
$\ol{X}$ et $\ol{Y}$ sont quotients de $X$ et $Y$, on voit que les
$B$-ensembles $\ol{X}$ et $\ol{Y}$ sont de la forme: $B/H\cprod \FF
\cprod \{\cdot\}$, avec $H$ contenant $\FF^*$ d'indice~$2$. Un tel $H$
est unique \`a conjugaison pr\`es dans $B$, donc $\ol{X}\cong\ol{Y}$
en tant que $B$-ensembles.

Comme $T'$ est cyclique, deux $T'$-ensembles transitifs sont
isomorphes si et seulement si ils ont m\^eme cardinal.  On a
$\#\ol{X}=\#\ol{Y}$. Alors, pour voir que $\ol{X}$ et $\ol{Y}$ sont
$T'$-isomorphes, il suffit de montrer qu'il sont tous les deux
$T'/\FF^*$-libres, \`a un point et \`a au plus une autre orbite
pr\`es.  Supposons qu'un $t$ non trivial dans $T'/\FF^*$ fixe un
\'el\'ement $\{P,Q\}$ de $\Sym^2(\PP^1(\FF))-\Delta$. Alors $t^2$ fixe
quatre points dans $\PP^1(\FFq)$, donc est scalaire.  Il existe au
plus un tel $t$ et on a $Q=tP$.  Comme $T'$ agit transitivement sur
$\PP^1(\FF)$, cela fait au plus une orbite non $T'/\FF^*$-libre dans
$\Sym^2(\PP^1(\FF))-\Delta$.  Notons $0'$ et $\infty'$ les deux points
fixes de $T'$ dans~$\PP^1(\FFq)$.  L'action de $T'$ sur la
$\FFq$-droite $0'$ dans $\FFq^2$ est un isomorphisme de $T'$
vers~$\FFq^*$; d\`es maintenant, nous verrons les $T'$-ensembles comme
des $\FFq^*$-ensembles via cet isomorphisme. L'action de $\FFq^*$ sur
la droite $\infty'$ est donn\'ee par~$\sigma$. Soit $e_1$ un
\'el\'ement non nul de $0'$, et $e_2:=\sigma(e_1)$ son image
dans~$\infty'$. Cette base de $\FFq^2$ donne une bijection entre
$\PP^1(\FFq)-\{0',\infty'\}$ et~$\FFq^*$, qui envoie la droite
$\FFq(e_1+ae_2)$ \`a $a\in\FFq^*$.  Via cette bijection, $\sigma$ sur
$\PP^1(\FFq)-\{0',\infty'\}$ correspond \`a $\phi\colon z\mapsto
1/\sigma(z)$ sur~$\FFq^*$. De m\^eme, pour $t$ dans $\FFq^*$, l'action
de $t$ devient $z\mapsto (\sigma(t)/t)z$.  Les orbites de $T'$ sur
$\FFq^*$ sont les fibres de la norme $\FFq^*\to\FF^*$, tandis que
$\phi$ agit sur l'ensemble de ces fibres par inversion de la
norme. L'action de $\phi$ sur l'orbite des \'el\'ements de norme $1$
est triviale, donc donne une orbite $T'/\FF^*$-libre dans $\ol{Y}$,
ainsi que les orbites de normes autres que $1$ et~$-1$.
\end{proof}

\begin{rem}
Nous avons formul\'e le Th\'eor\`eme~\ref{thm2.1} pour des objets avec
action de $G$ dans des cat\'egories additives $\QQ$-lin\'eaires
admettant des invariants, pour pouvoir l'appliquer \'egalement aux
motifs de Chow, par exemple ceux associ\'es aux espaces de formes
modulaires de poids au moins deux, construits comme
dans~\cite{Scholl2}, ou dans d'autres cat\'egories $\QQ$-lin\'eaires
pseudo-ab\'eliennes.  Rappelons (par exemple \cite[\S1]{Scholl1})
qu'une cat\'egorie additive est dite pseudo-ab\'elienne si pour tout
objet $M$ tout idempotent dans $\End(M)$ admet un noyau (ou, ce qui
est \'equivalent, une image). Pour un groupe fini $H$ op\'erant sur un
objet $M$ l'objet $M^H$ est l'image de l'idempotent $(1/\# H)\sum_h
h$.
\end{rem}
\begin{rem}
L'id\'ee d'utiliser la th\'eorie des repr\'esentations d'un groupe
fini $G$ pour en d\'eduire des iso\-g\'e\-nies entre vari\'et\'es
ab\'eliennes n'est certainement pas nouvelle; voir \cite{KaRo1}
et~\cite{KaRo2}.  La preuve du Th\'eor\`eme~\ref{thm2.1} montre que ce
genre de r\'esultats reste vrai dans le cadre plus g\'en\'eral des
cat\'ego\-ries additives $\QQ$-lin\'eaires pseudo-ab\'eliennes.  Dans
\cite[\S5]{KaRo1} on trouve une relation pour certains sous-groupes de
$\SL_2(\FF_p)/\{1,-1\}$, mais pas les relations du
Th\'eor\`eme~\ref{thm2.1}. L'int\'er\^et du r\'esultat de Chen est que
l'on obtient des renseignements sur la jacobienne de $X(p)_\nonsplit$
en termes d'objets d\'ej\`a mieux compris.
\end{rem}
\begin{rem}
Lo\"\i c Merel a pos\'e la question de savoir si la correspondance
constitu\'ee des deux morphismes quotients $X\to X/N'$ et $X\to X/N$
induit un morphisme de la jacobienne de $X/N'$ vers celle de $X/N$
dont le noyau est fini. La r\'eponse d\'epend de la position relative
des sous-groupes $N$ et $N'$ dans~$G$. En effet, \`a conjugaison
pr\`es, il existe un unique couple $(T,T')$ avec $T$ d\'eploy\'e et
$T'$ non d\'eploy\'e, telle que $N\cap N'$ soit de cardinal $4(p-1)$.
Pour cette configuration, Chen a montr\'e que le noyau en question est
fini; voir~\cite{Chen2}.
\end{rem}

\noindent
Nous nous int\'eressons maintenant \`a la question de savoir quels
d\'enominateurs sont essentiels dans le Th\'eor\`eme~2. Plus
pr\'ecis\'ement, on veut remplacer la condition ``$\QQ$-lin\'eaire''
sur la cat\'ego\-rie $\calC$ par ``$\ZZl$-lin\'eaire'' avec $l$
premier et $\ZZl$ le localis\'e de $\ZZ$ en~$l$. La th\'eorie de
l'induction de Conlon \cite[\S81B]{CR} nous permet de donner une
r\'eponse compl\`ete.

\begin{theorem}
Soit $q$ le cardinal de~$\FF$.  Soit $l$ un nombre premier. Si $l$ ne
divise pas $q^2-1$, alors le Th\'eor\`eme~2 reste vrai avec $\QQ$
remplac\'e par~$\ZZl$.  Si $l$ ne divise pas $q-1$, alors la partie
(2) du Th\'eor\`eme~2 reste vraie avec $\QQ$ remplac\'e par~$\ZZl$.
\end{theorem}
\begin{proof}
Soit $l$ premier. Un groupe fini $C$ est dit cyclique modulo $l$ si
$C$ admet un quotient cyclique de noyau un $l$-groupe. Pour $H$ un
groupe fini et $U$, $V$ deux $H$-ensembles finis la th\'eorie de
l'induction de Conlon implique que les $H$-modules $\ZZl[U]$ et
$\ZZl[V]$ sont isomorphes si et seulement si pour tout sous-groupe $C$
cyclique modulo $l$ de $H$, les $C$-ensembles $U$ et $V$ sont
isomorphes.  Pour voir ``seulement si" on applique (81.28) dans
\cite{CR} avec $R=\FF_l$, et pour ``si" on utilise (81.25) avec
$R=\FF_l$, et les arguments dans (81.17) avec $R=\ZZl$.

Supposons que $l$ ne divise pas l'ordre de~$G$. Alors les sous-groupes
de $G$ cycliques modulo $l$ sont simplement les sous-groupes
cycliques.  Nous avons d\'ej\`a vu que $X$ et $Y$, ainsi que $\ol{X}$
et $\ol{Y}$, sont $C$-isomorphes pour tout sous-groupe cyclique $C$
de~$G$.


Consid\'erons le cas o\`u $l$ est la caract\'eristique de~$\FF$.
Comme chaque \'el\'ement dans $G$ d'ordre $l$ a une point fixe unique
dans $\PP^1(\FF)$, tout sous-groupe $H$ de $G$ cyclique modulo $l$ est
soit cyclique soit un conjugu\'e d'un sous-groupe de~$B$.  Lors de la
preuve du Th\'eor\`eme~2 nous avons vu que $X$ et $Y$, ainsi que
$\ol{X}$ et $\ol{Y}$, sont des $B$-ensembles isomorphes. Ceci termine
la preuve pour ce~$l$.

La derni\`ere chose \`a montrer est que
$\ZZl[\ol{X}]\cong\ZZl[\ol{Y}]$ en tant que $G$-modules pour les
$l\neq2$ divisant $q+1$. Supposons $l$ comme cela.  Notons que
$N'/\FF^*$ est un groupe di\'edrale d'ordre $2(q+1)$.  Les
$N'/\FF^*$-ensembles $\ol{X}$ et $\ol{Y}$ ont chacun un point fixe
unique; notons $X'$ et $Y'$ les compl\'ements. Soit $S$ un
$2$-sous-groupe de Sylow $N'/\FF^*$, et $H$ le sous-groupe maximal
d'ordre impair de~$N'/\FF^*$; ainsi $N'/\FF^*$ est le produit
semi-direct de $H$ par~$S$. Comme $\QQ[X']\cong_{N'}\QQ[Y']$ on a
$\QQ[H\backslash X']\cong_{S}\QQ[H\backslash Y']$. Comme $l$ ne divise
pas $\# S$, on a $\ZZl[H\backslash X']\cong_{S}\ZZl[H\backslash Y']$
par Conlon.  Nous avons d\'ej\`a vu que les stabilisateurs dans
$T'/\FF^*$ des \'el\'ements de $X'$ et de $Y'$ sont d'ordre un ou
deux. Donc les stabilisateurs dans $N'/\FF^*$ ont des conjugu\'es
dans~$S$. Il en r\'esulte que $X'$ et $Y'$ sont les induits des
$S$-ensembles $H\backslash X'$ et $H\backslash Y'$
(cf. \cite[\S80B]{CR}).  Par induction de modules, $\ZZl[X']$ et
$\ZZl[Y']$ sont isomorphes en tant que $N'$-modules. Par Conlon, $X'$
et $Y'$ sont isomorphes en tant que $C$-ensembles pour tout
sous-groupe $C$ de $N'$ cyclique modulo $l$. Mais tout sous-groupe
cyclique modulo $l$ de $G$ est soit cyclique soit un conjugu\'e d'un
sous-groupe de~$N'$. Encore par Conlon,
$\ZZl[\ol{X}]\cong_G\ZZl[\ol{Y}]$.
\end{proof}
\begin{rem}
Soit $l$ un nombre premier qui divise $q^2-1$.  Il existe un
sous-groupe $D$ de $N/\FF^*$ ou de $N'/\FF^*$, di\'edral d'ordre $2l$.
Alors $D$ fixe plus de points dans $X$ que dans $Y$, et, si $l$ divise
$q-1$, plus dans $\ol{X}$ que dans $\ol{Y}$.  Par Conlon, on voit que
pour la cat\'egorie des $\ZZl$-modules les parties (1) et (2) du
Th\'eor\`eme 2 sont vraies seulement si $l$ satisfait les conditions
du Th\'eor\`eme~3.
\end{rem}

\begin{theorem}
Notations comme dans le Th\'eor\`eme~1. Soit $l$ un nombre premier.
Si $l$ ne divise pas $p^2-1$, alors il existe une isog\'enie
$\jac(X/T)\to\jac(X/T')\times\jac(X/B)^2$ de degr\'e premier \`a~$l$.
Si $l$ ne divise pas $p-1$, alors il existe une isog\'enie
$\jac(X/N)\to\jac(X/N')\times \jac(X/B)$ de degr\'e premier \`a~$l$.
\end{theorem}
\begin{proof}
On note d'abord que pour des objets $A$ et $B$ de $\calA$ on a
$\ZZl\otimes A\cong\ZZl\otimes B$ si et seulement s'il existe une
isog\'enie $A\to B$ de degr\'e premier \`a~$l$. Sous les hypoth\`eses
faites sur $l$, le Th\'eor\`eme~3 donne l'existence d'isomorphismes
$(\ZZl\otimes\jac(X))^T\to(\ZZl\otimes\jac(X))^{T'}
\oplus(\ZZl\otimes\jac(X))^{B,2}$ et
$(\ZZl\otimes\jac(X))^N\to(\ZZl\otimes\jac(X))^{N'}
\oplus(\ZZl\otimes\jac(X))^B$. Pour terminer, il suffit de voir que
les noyaux des morphismes $\jac(X/H)\to\jac(X)$ sont d'ordre premier
aux $l$ concern\'es, pour $H$ parmi $T$, $T'$, $B$, $N$ et~$N'$. Mais
ce noyau est le dual de Cartier du groupe du plus grand
sous-rev\^etement \'etale ab\'elien de
$(X/H)_\Qbar=Y/(H\cap\SL_2(\FF_p))$ dans $Y\to (X/H)_\Qbar$, o\`u $Y$
est une composante connexe de $X_\Qbar$ (sur $\CC$ cela se voit \`a
l'aide des groupes fondamentaux). Pour $p<5$ on a $\jac(X)=0$, nous
supposons donc que $p\geq5$.  Pour $N$ et $N'$, les noyaux en question
sont alors des $2$-groupes. Pour $T$ et $B$ on obtient des quotients
de $\FF_p^*$, et $T'$ donne un groupe cyclique d'ordre divisant $p+1$.
\end{proof}

\smallskip\noindent {\bf Remerciements.} Nous tenons \`a remercier
Laurent Moret-Bailly pour une lecture critique de ce texte.

\vfill
\bigskip
\hbox to \hsize{
  \vbox{
    \hbox{Bart de Smit\hss}
    \hbox{Mathematisch Instituut, Universiteit Leiden\hss}
    \hbox{Postbus 9512, 2300~RA Leiden, Pays-Bas\hss}
    \smallskip
    \hbox{\tt desmit@math.leidenuniv.nl\hss}
    \hbox{tel.: +31(0)715277144\hss}
    \hbox{fax: +31(0)715277101\hss}
   }
\hss
  \vbox{
    \hbox{Bas Edixhoven \hss}
    \hbox{IRMAR, Universit\'e de Rennes 1\hss}
    \hbox{35042 Rennes Cedex, France\hss}
    \smallskip
    \hbox{\tt edix@maths.univ-rennes1.fr}
    \hbox{tel.: +33(0)299286018\hss}
    \hbox{fax: +33(0)299286790\hss}
   }
}

\end{document}